\documentclass[12pt]{article}%
\usepackage{amsfonts}
\usepackage{amsmath,amssymb}
\usepackage[applemac]{inputenc}
\usepackage{amsmath,amssymb,fullpage}
\usepackage{color}
\usepackage{amsmath}
\usepackage{amssymb}
\usepackage{graphicx}
\usepackage{tikz}
\providecommand{\U}[1]{\protect\rule{.1in}{.1in}}
\usetikzlibrary{arrows,shapes,automata,petri}

\newcommand{\dproof}{\noindent {Proof.} \quad}
\newcommand{\fproof}{\hfill $\square$ \bigskip}
\newtheorem{definition}{Definition}[section]

\newtheorem{theorem}[definition]{Theorem}
\newtheorem{problem}[definition]{Problem}
\newtheorem{remark}[definition]{ \it Remark}

\newtheorem{lemma}[definition]{Lemma}

\numberwithin{equation}{section}

\def\1B{\text{1\!\!I}}

\begin{document}

\title{Singular control of SPDEs with space-mean dynamics}
\author{Nacira AGRAM$^{1}$ Astrid HILBERT$^{1}$ and Bernt ØKSENDAL$^{2}$}
\date{3 May 2019}
\maketitle

\begin{abstract}
We consider the problem of optimal singular control of a stochastic partial
differential equation (SPDE) with space-mean dependence. Such systems are
proposed as models for population growth in a random environment. We obtain
sufficient and necessary maximum principles for such control problems. The
corresponding adjoint equation is a reflected backward stochastic partial
differential equation (BSPDE) with space-mean dependence. We prove existence
and uniqueness results for such equations. As an application we study optimal
harvesting from a population modelled as an SPDE with space-mean dependence.

\end{abstract}



\footnotetext[1]{Department of Mathematics, Linnaeus University (LNU),
Sweden.\newline Emails: \texttt{nacira.agram@lnu.se, astrid.hilbert@lnu.se.}}

\footnotetext[2]{Department of Mathematics, University of Oslo, Norway.
\newline Email: \texttt{oksendal@math.uio.no.}
\par
This research was carried out with support of the Norwegian Research Council,
within the research project Challenges in Stochastic Control, Information and
Applications (STOCONINF), project number 250768/F20.}

\section{Introduction}

We start by a motivation for the problem that will be studied in this
paper:\newline Consider a problem of optimal harvesting from a fish population
in a lake $D$. We assume that the density $u(t,x)$ of the population at time
$t\in\lbrack0,T]$ and at the point $x\in D$ is modelled by a \emph{stochastic
reaction-diffusion equation with neighbouring interactions}. By this we mean a
stochastic partial differential equation of the form
\[
\left\{
\begin{array}
[c]{l}%
du(t,x)=\left[  \dfrac{1}{2}\Delta u(t,x)+\alpha\overline{u}(t,x)\right]
dt+\beta\overline{u}(t,x)dB(t)-\lambda_{0}\xi(dt,x);\quad(t,x)\in(0,T)\times
D\\
\text{ }u(0,x)=u_{0}(x)>0;\quad x\in D,\\
\text{ }u(t,x)=u_{1}(t,x)\geq0;\quad(t,x)\in(0,T)\times\partial D,
\end{array}
\right.
\]
where $\overline{u}(t,x)$ is the space-averaging operator
\begin{equation}
\overline{u}(t,x)=\frac{1}{V(K_{\theta})}\int_{K_{\theta}}u(x+y)dy, \label{G}%
\end{equation}
where $V(\cdot)$ denotes Lebesgue volume and
\[
K_{\theta}=\{y\in\mathbb{R}^{n};|y|<\theta\}
\]
is the ball of radius $r>0$ in $\mathbb{R}^{d}$ centered at $0$,where $D$ is a
bounded Lipschitz domain in $\mathbb{R}^{d}$ and $u_{0}(x),u_{1}(t,x)$ are
given deterministic functions. \newline In the above $B(t)=B(t,\omega
);(t,\omega)\in\lbrack0,\infty)\times\Omega,$ is an $m$-dimensional Brownian
motion on a filtered probability space $(\Omega,\mathbb{F}=\{\mathcal{F}%
_{t}\}_{t\in\lbrack0,\infty)},\mathbb{P})$. Moreover, $\alpha$, $\beta$ and
$\lambda_{0}>0$ are given constants and
\[
\Delta=\sum\limits_{i=1}^{d}{\frac{\partial^{2}}{{\partial x_{i}^{2}}}}%
\]
is the Laplacian differential operator on $\mathbb{R}^{d}$.\newline We may
regard $\xi(dt,x)$ as the harvesting effort rate, and $\lambda_{0}>0$ as the
harvesting efficiency coefficient. The performance functional is assumed to be
of the form
\begin{equation}
J(\xi)=\mathbb{E}\left[  \int_{D}\int_{0}^{T}(h_{0}(t,x)u(t,x)-c(t,x))\xi
(dt,x)dx+\int_{D}h_{0}(T,x)u(T,x)dx\right]  , \label{eq248}%
\end{equation}
where $h_{0}(t,x)>0$ is the unit price of the fish and $c(t,x)$ is the unit
cost of energy used in the harvesting and $T>0$ is a fixed terminal time. Thus
$J(\xi)$ represents the expected total net income from the harvesting. The
problem is to maximise $J(\xi)$ over all (admissible) harvesting strategies
$\xi(t,x)$.

\begin{remark}
This population growth model, which was first introduced in Agram et al
\cite{AHO}, is a generalisation of the classical stochastic reaction-diffusion
model, in that we have added the term $\overline{u}(t,x)$ which represents an
average of the neighbouring densities. Thus our model allows for the growth at
a point to depend on interactions from the whole vicinity. This space-mean
interaction is different from the pointwise interaction represented by the Laplacian.
\end{remark}

\noindent The problem above turns out to be related to a problem of the
following form:\newline Let $\phi(x)=\phi(x,\omega)$ be an $\mathcal{F}_{T}%
$-measurable $H=L^{2}(D)$-valued random variable. Let
\[
g:[0,T]\times D\times\mathbb{R}\times\mathbb{R}^{m}\rightarrow\mathbb{R}%
\]
be a given measurable mapping and $L(t,x):[0,T]\times D\rightarrow\mathbb{R}$
a given continuous function. Consider the problem to find an $\mathbb{F}%
$-adapted random fields $Y(t,x)\in\mathbb{R},Z(t,x)\in\mathbb{R}^{m}%
,\xi(t,x)\in\mathbb{R}^{+}$ left-continuous and increasing with respect to
$t$, such that
\begin{equation}
\left\{
\begin{array}
[c]{l}%
dY(t,x)=-AY(t,x)dt-F(t,x,Y(t,x),Y(t,\cdot),Z(t,x))dt+Z(t,x)dB(t)\\
\quad\quad\quad\quad\quad\quad\quad\quad\quad\quad-\xi(dt,x),t\in(0,T),\\
Y(t,x)\geq L(t,x),\\
\int_{0}^{T}\int_{D}(Y(t,x)-L(t,x))\xi(dt,x)dx=0,\\
Y(T,x)=\phi(x)\quad\text{a.s.,}%
\end{array}
\right.  \label{a}%
\end{equation}
where $A$ is a second order linear partial differential operator. We call the
equation \eqref{a} a \emph{reflected stochastic partial differential equation
(SPDE) with space-mean dynamics}. We will come back to this equation in the
last section.

\section{The optimization problem}

We now give a general formulation of the problem discussed in the
Introduction:\newline

\noindent Let $T>0$ and let $D\subset%
%TCIMACRO{\U{211d} }%
%BeginExpansion
\mathbb{R}
%EndExpansion
^{n}$ be an open set with $C^{1}$ boundary $\partial D.$ Specifically, we
assume that the state $u(t,x)$ at time $t\in\lbrack0,T]$ and at the point
$x\in\overline{D}:=D\cup\partial D$ satisfies%
\begin{equation}%
\begin{cases}
du(t,x) & =Au(t,x)dt+b(t,x,u(t,x),u(t,\cdot))dt+\sigma(t,x,u(t,x),u(t,\cdot
))dB(t)\\
& \text{ \ \ \ \ }+f(t,x,u)\xi(dt,x);\quad(t,x)\in(0,T)\times D,\\
u(0^{-},x) & =u_0(x);\quad x\in\overline{D},\\
u(t,x) & =u_1(t,x);\quad(t,x)\in(0,T)\times\partial D.
\end{cases}
\label{eq2.1}
\end{equation}
Here $B=\{B(t)\}_{t\in\lbrack0,T]}$ is a $d$-dimensional Brownian motion,
defined in a complete filtered probability space $(\Omega,\mathcal{F}%
,\mathbb{F},\mathbb{P}).$ The filtration $\mathbb{F=}\left\{  \mathcal{F}%
_{t}\right\}  _{t\geq0}$ is assumed to be the $\mathbb{P}$-augmented
filtration generated by $B$.\newline We denote by $A$ the second order partial
differential operator acting on $x$ given by%
\begin{equation}
A\phi(x)=\sum_{i,j=1}^{n}\alpha_{ij}(x)\frac{\partial^{2}\phi}{\partial
x_{i}\partial x_{j}}+\sum_{i=1}^{n}\beta_{i}(x)\frac{\partial\phi}{\partial
x_{i}};\quad\phi\in C^{2}(%
%TCIMACRO{\U{211d} }%
%BeginExpansion
\mathbb{R}
%EndExpansion
^{n}), \label{opr}%
\end{equation}
where $(\alpha_{ij}(x))_{1\leq i,j\leq n}$ is a given nonnegative definite
$n\times n$ matrix with entries $\alpha_{ij}(x)\in C^{2}(D)\cap C(\overline
{D})$ for all $i,j=1,2,...,n$ and $\beta_{i}(x)\in C^{2}(D)\cap C(\overline
{D})$ for all $i=1,2,...,n.$\newline Let $L(\mathbb{R}^{n})$ denote the set of
real measurable functions on $\mathbb{R}^{n}$. For each $t,x,u,\zeta$ the
functions
\begin{align*}
\varphi &  \mapsto b(t,x,u,\varphi):[0,T]\times D\times%
%TCIMACRO{\U{211d} }%
%BeginExpansion
\mathbb{R}
%EndExpansion
\times L(\mathbb{R}^{n})\rightarrow%
%TCIMACRO{\U{211d} }%
%BeginExpansion
\mathbb{R}
%EndExpansion
,\\
\varphi &  \mapsto\sigma(t,x,u,\varphi):[0,T]\times D\times%
%TCIMACRO{\U{211d} }%
%BeginExpansion
\mathbb{R}
%EndExpansion
\times L(\mathbb{R}^{n})\rightarrow%
%TCIMACRO{\U{211d} }%
%BeginExpansion
\mathbb{R},\\
u &  \mapsto f(t,x,u):[0,T]\times D\times%
%TCIMACRO{\U{211d} }%
%BeginExpansion
\mathbb{R}
%EndExpansion
\rightarrow%
%TCIMACRO{\U{211d} }%
%BeginExpansion
\mathbb{R},
\end{align*}
are $C^{1}$ functionals on $L^{2}(D)=L^{2}(D,m)$, where $dm(x)=dx$ is the
Lebesgue measure on $\mathbb{R}^{n}$. Here $Au(t,x)$ is interpreted in the
sense of distribution. Thus $u$ is understood as a weak (mild) solution to
(\ref{eq2.1}), in the sense that
\begin{align}
u(0,x)&=u_{0}(x)+\int_{0}^{t}P_{s}^{A}b(s,\cdot,u(s,x),u(s,\cdot))(x)ds+\int
_{0}^{t}P_{s}^{A}\sigma(s,\cdot,u(s,x),u(s,\cdot))(x)dB(s)\nonumber\\
&+\int_{0}^{t}P_{s}^{A}f(s,\cdot,u(s,x))(x)\xi(ds,x),
\end{align}
where $P_{t}^{A}$ is the semigroup associated to the operator $A$. Thus
we see that we can in the usual way apply the Itô formula to such
SPDEs.\newline Moreover, the adjoint operator $A^{\ast}$ of an operator $A$ on
$C_{0}^{\infty}(\mathbb{R})$ is defined by the identity
\[
(A\phi,\psi)=(\phi,A^{\ast}\psi),\,\,\,\text{ for all }\phi,\psi\in
C_{0}^{\infty}(\mathbb{R}),
\]
where $\langle\phi_{1},\phi_{2}\rangle_{L^{2}(\mathbb{R})}:=(\phi_{1},\phi
_{2})=\int_{\mathbb{R}}\phi_{1}(x)\phi_{2}(x)dx$ is the inner product in
$L^{2}(\mathbb{R}).$ In our case we have
\[
A_{x}^{\ast}\phi(x)=\sum_{i,j=1}^{n}\frac{\partial^{2}}{\partial x_{i}\partial
x_{j}}(\alpha_{ij}(x)\phi(x))-\sum_{i=1}^{n}\frac{\partial}{\partial x_{i}%
}(\beta_{i}(x)\phi(x));\quad\phi\in C^{2}(%
%TCIMACRO{\U{211d} }%
%BeginExpansion
\mathbb{R}
%EndExpansion
^{n}).
\]
We interpret $u$ as a weak (variational) solution to \eqref{eq2.1}, in the
sense that for $\phi\in C_{0}^{\infty}(D),$
\begin{align*}
\langle u(t),\phi\rangle_{L^{2}(D)}  &  =\langle\eta(\cdot),\phi\rangle
_{L^{2}(D)}+\int_{0}^{t}\langle u(s),A_{x}^{\ast}\phi\rangle ds+\int_{0}^{t}\langle b(s,u(s)),\phi\rangle_{L^{2}(D)}ds\\
&+\int_{0}^{t}\langle\sigma(s,u(s)),\phi\rangle_{L^{2}(D)}dB(s)+\int_{0}^{t}\langle f(s.u(s)),\phi\rangle_{L^{2}(D)} \xi(ds,x),
\end{align*}
where $\langle\cdot,\cdot\rangle$ represents the duality product between
$W^{1,2}(D)$ and $W^{1,2}(D)^{\ast}$, with $W^{1,2}(D)$ the Sobolev space of
order $1$. In the above equation, we have not written all the arguments of
$b,\sigma,\gamma$, for simplicity. \newline We want to maximize the
performance functional $J(\xi),$ given by
\begin{equation}%
\begin{array}
[c]{c}%
J(\xi)=\mathbb{E}\Big[%
%TCIMACRO{\dint _{0}^{T}}%
%BeginExpansion
{\displaystyle\int_{0}^{T}}
%EndExpansion%
%TCIMACRO{\dint _{D}}%
%BeginExpansion
{\displaystyle\int_{D}}
%EndExpansion
h_{0}(t,x,u(t,x),u(t,\cdot))dxdt+%
%TCIMACRO{\dint _{0}^{T}}%
%BeginExpansion
{\displaystyle\int_{0}^{T}}
%EndExpansion%
%TCIMACRO{\dint _{D}}%
%BeginExpansion
{\displaystyle\int_{D}}
%EndExpansion
h_{1}(t,x,u(t,x),u(t,\cdot))\xi(dt,dx)\\
+%
%TCIMACRO{\dint _{D}}%
%BeginExpansion
{\displaystyle\int_{D}}
%EndExpansion
g(x,u(T,x),u(T,\cdot))dx\Big],
\end{array}
\label{ju}%
\end{equation}
over all $\xi\in\mathcal{A}$, where $\mathcal{A}$ is the set of all adapted
processes $\xi(t,x)$ that are nondecreasing and left continuous with respect
to $t$ for all $x$, with $\xi(0,x)=0,$ $\xi(T,x)<\infty$ and such that
$J(\xi)<\infty.$ We call $\mathcal{A}$ the set of admissible singular
controls. Thus we want to find $\widehat{\xi}\in\mathcal{A},$ such that
\begin{equation}
J(\widehat{\xi})=\sup_{\xi\in\mathcal{A}}J(\xi). \label{eq2.4}%
\end{equation}
For each $t,x,u$ we assume that the functions $\varphi\mapsto h_{0}%
(t,x,u,\varphi):[0,T]\times D\times%
%TCIMACRO{\U{211d} }%
%BeginExpansion
\mathbb{R}
%EndExpansion
\times L(\mathbb{R}^{n})\rightarrow%
%TCIMACRO{\U{211d} }%
%BeginExpansion
\mathbb{R}
%EndExpansion
,$ and $\varphi\mapsto g(x,u,\varphi):D\times%
%TCIMACRO{\U{211d} }%
%BeginExpansion
\mathbb{R}
%EndExpansion
\times L(\mathbb{R}^{n})\rightarrow%
%TCIMACRO{\U{211d} }%
%BeginExpansion
\mathbb{R}
%EndExpansion
,$ are $C^{1}$ functionals on $L^{2}(D)$.\newline The Hamiltonian $H$ is
defined by%

\begin{align}
\label{Ham} &  H(t,x,u,\varphi,p,q)(dt,\xi(dt,x))= H_{0}(t,x,u,\varphi,p,q)dt
+ H_{1}(t,x,u,\varphi,p) \xi(dt,x).
\end{align}
where
\begin{align}
H_{0}(t,x,u,\varphi,p,q)=  &  = h_{0}(t,x,u,\varphi)+b(t,x,u,\varphi
)p+\sigma(t,x,u,\varphi)q \label{H_0}%
\end{align}
and
\begin{align}
H_{1}(t,x,u,\varphi,p)=  &  = f(t,x)p+ h_{1}(t,x,u,\varphi) \label{H_1}%
\end{align}
We assume that $H,f,b,\sigma,\gamma$ and $g$ admit Fréchet derivatives with
respect to $u$ and $\varphi.$ \newline In general, if $h:L^{2}(D)\mapsto
L^{2}(D)$ is Fréchet differentiable, we denote its Fréchet derivative
(gradient) at $\varphi\in L^{2}(D)$ by $\nabla_{\varphi}h$, and we denote the
action of $\nabla_{\varphi}h$ on a function $\psi\in L^{2}(D)$ by
$\left\langle \nabla_{\varphi}h,\psi\right\rangle $.

\begin{definition}
We say that the Fréchet derivative $\nabla_{\varphi}h$ of a map $h:L^{2}%
(D)\mapsto L^{2}(D)$ has a \emph{dual function} $\nabla_{\varphi}^{\ast}h\in
L^{2}(D\times D)$ if
\begin{equation}
\left\langle \nabla_{\varphi}h,\psi\right\rangle (x)=\int_{D}\nabla_{\varphi
}^{\ast}h(x,y)\psi(y)dy;\quad\text{ for all }\psi\in L^{2}(D). \label{eq2.6a}%
\end{equation}

\end{definition}

\noindent By Fubini's theorem, we get%
\begin{equation}
\overline{\nabla}_{\varphi}^{\ast}h(x)=\int_{D}\nabla_{\varphi}^{\ast
}h(y,x)dy. \label{eq2.8a}%
\end{equation}
We associate to the Hamiltonian the following reflected BSPDE
\begin{equation}
\left\{
\begin{array}
[c]{ll}%
dp(t,x) & =-A_{x}^{\ast}p(t,x)dt-\left\{  \frac{\partial H_{0}}{\partial
u}(t,x)+\overline{\nabla}_{\varphi}^{\ast}H_{0}(t,x)\right\}  dt\\
& -\left\{  \frac{\partial H_{1}}{\partial u}(t,x)+\overline{\nabla}_{\varphi
}^{\ast}H_{1}(t,x)\right\}  \xi(dt,x)\\
& \text{ \ }+q(t,x)dB(t);\quad(t,x)\in(0,T)\times D,\\
p(T,x) & =\frac{\partial g}{\partial u}(T,x)+\overline{\nabla}_{\varphi}%
^{\ast}g(T,x);\quad x\in D,\\
p(t,x) & =0;\quad(t,x)\in(0,T)\times\partial D,
\end{array}
\right.  \label{p}%
\end{equation}
where we have used the simplified notation
\[
H_{i}(t,x)=H_{i}(t,x,u,\varphi,p,q)|_{u=u(t,x),\varphi=u(t,\cdot
),p=p(t,x),q=q(t,x)},\quad i=0,1
\]
and similarly with $g$.

\subsection{A sufficient maximum principle}

We now formulate a sufficient version ( a verification theorem) of the maximum
principle for the optimal control of the problem \eqref{eq2.1}-\eqref{eq2.4}.

\begin{theorem}
[Sufficient Maximum Principle]Suppose $\widehat{\xi}\in\mathcal{A}$, with
corresponding \newline$\widehat{u}(t,x),\widehat{p}(t,x),\widehat{q}(t,x).$
Suppose the functions $(u,\varphi)\mapsto g(x,u,\varphi)$ and\newline%
$(u,\varphi,\xi)\mapsto H(t,x,u,\varphi,\widehat{p}(t,x),\widehat
{q}(t,x))(dt,\xi(dt,dx))$ are concave for each $(t,x)\in(0,T)\times D$.
Moreover, suppose that
\[
\widehat{\xi}(dt,x)\in\arg\max_{\widehat{\xi}\in\mathcal{A}}H(t,x,\widehat
{u}(t,x),\widehat{u}(t,\cdot),\widehat{p}(t,x),\widehat{q}(t,x))(dt,\xi
(dt,x));
\]
i.e.,%
\begin{align}
&  \left\{  H(t,x,\widehat{u}(t,x),\widehat{u}(t,\cdot),\widehat
{p}(t,x),\widehat{q}(t,x)\right\}  \xi(dt,x)\label{max-cond}\\
&  \leq\left\{  H(t,x,\widehat{u}(t,x),\widehat{u}(t,\cdot),\widehat
{p}(t,x),\widehat{q}(t,x)\right\}  \widehat{\xi}(dt,x);\text{ for all }\xi
\in\mathcal{A}.\nonumber
\end{align}

Then $\widehat{\xi}$ is an optimal singular control.
\end{theorem}

\dproof
Consider
\[
J(\xi)-J(\widehat{\xi})=I_{1}+I_{2}+I_{3},
\]
where
\[
I_{1}=\mathbb{E}\left[  \int_{0}^{T}{\int_{D}}\{h_{0}(t,x,u(t,x),u(t,\cdot
))-h_{0}(t,x,\widehat{u}(t,x),\widehat{u}(t,\cdot))\}dxdt\right]  ,
\]%
\[
I_{2}=\mathbb{E}\left[  \int_{0}^{T}{\int_{D}}\{h_{1}(t,x,u(t,x),u(t,\cdot
))\xi(dt,x)-\int_{0}^{T}{\int_{D}}h_{1}(t,x,\widehat{u}(t,x),\widehat
{u}(t,\cdot))\widehat{\xi}(dt,x)\right]  ,
\]
and
\begin{equation}
I_{3}={\int_{D}}\mathbb{E}\left[  g(x,u(T,x),u(T,\cdot))-g(x,\widehat
{u}(T,x),\widehat{u}(T,\cdot))\right]  dx. \label{I3}%
\end{equation}
By concavity on $g$ together with the identity \eqref{eq2.6a}-\eqref{eq2.8a},
we get
\begin{align*}
I_{3}  &  \leq{\int_{D}}\mathbb{E}\left[  \frac{\partial\widehat{g}}{\partial
u}(T,x)(u(T,x)-\widehat{u}(T,x))+\left\langle \nabla_{\varphi}\widehat
{g}(T,x),u(T,\cdot)-\widehat{u}(T,\cdot)\right\rangle \right]  dx\\
&  ={\int_{D}}\mathbb{E}\left[  \frac{\partial\widehat{g}}{\partial
u}(T,x)(u(T,x)-\widehat{u}(T,x))+\overline{\nabla}_{\varphi}^{\ast}\widehat
{g}(T,x)(u(T,x)-\widehat{u}(T,x))\right]  dx\\
&  ={\int_{D}}\mathbb{E}\left[  \widehat{p}(T,x)(u(T,x)-\widehat
{u}(T,x))\right]  dx\\
&  ={\int_{D}}\mathbb{E}\left[  \widehat{p}(T,x)\widetilde{u}(T,x)\right]  dx,
\end{align*}
where $\tilde{u}(t,x)=u(t,x)-\hat{u}(t,x);t\in\lbrack0,T]$.

\noindent Applying the Itô formula to $\widehat{p}(t,x)\widetilde{u}(t,x)$, we
have
\begin{align}
I_{3}  &  \leq\mathbb{E}\Big[{\int_{0}^{T}}{\int_{D}}\Big(\widehat
{p}(t,x)\Big\{A_{x}\widetilde{u}(t,x)+\widetilde{b}(t,x)\Big\}-\widetilde
{u}(t,x)A_{x}^{\ast}\widehat{p}(t,x)\nonumber\\
&  -\widetilde{u}(t,x)\Big\{\frac{\partial\widehat{H}_{0}}{\partial
u}(t,x)+\overline{\nabla}_{\varphi}^{\ast}\widehat{H}_{0}(t,x)\Big\}+\widehat
{q}(t,x)\widetilde{\sigma}(t,x)\Big)dxdt\nonumber\\
&  +{\int_{0}^{T}}{\int_{D}}\Big(\widehat{p}(t,x)f(t,x)\Big(\xi(dt,x)-\widehat
{\xi}(dt,x)\Big)-\widetilde{u}(t,x)\Big\{\frac{\partial\widehat{H}_{1}%
}{\partial u}(t,x)+\overline{\nabla}_{\varphi}^{\ast}\widehat{H}%
_{1}(t,x)\Big\}\Big)\xi(dt,x)\Big]. \label{1}%
\end{align}
By the first Green formula (see e.g. Wloka \cite{W}, page 258) there exist
first order boundary differential operators $A_{1},A_{2},$ such that%
\begin{align*}
&  \int_{D}\left\{  \widehat{p}(t,x)A_{x}\widetilde{u}(t,x)-\widetilde
{u}(t,x)A_{x}^{\ast}\widehat{p}(t,x)\right\}  dx\\
&  =\int_{\partial D}\left\{  \widehat{p}(t,x)A_{1}\widetilde{u}%
(t,x)-\widetilde{u}(t,x)A_{2}\widehat{p}(t,x)\right\}  dS,
\end{align*}
where the last integral is the surface integral over $\partial D$. We have
that
\begin{equation}
\widetilde{u}(t,x)=\widehat{p}(t,x)\equiv0, \label{**}%
\end{equation}
for all $(t,x)\in(0,T)\times\partial D.$\newline Substituting $\left(
\ref{**}\right)  $ in $\left(  \ref{1}\right)  $, yields%
\begin{align}
I_{3}  &  \leq\mathbb{E}\Big[{\int_{0}^{T}}{\int_{D}}\Big(\widehat
{p}(t,x)\widetilde{b}(t,x)-\widetilde{u}(t,x)\Big\{\frac{\partial\widehat
{H}_{0}}{\partial u}(t,x)+\overline{\nabla}_{\varphi}^{\ast}\widehat{H}%
_{0}(t,x)\Big\}+\widehat{q}(t,x)\widetilde{\sigma}(t,x)\Big)dxdt\nonumber\\
&  +{\int_{0}^{T}}{\int_{D}}\Big(\widehat{p}(t,x)f(t,x)\Big(\xi(dt,x)-\widehat
{\xi}(dt,x)\Big)-\widetilde{u}(t,x)\Big\{\frac{\partial\widehat{H}_{1}%
}{\partial u}(t,x)+\overline{\nabla}_{\varphi}^{\ast}\widehat{H}%
_{1}(t,x)\Big\}\xi(dt,x)\Big)dx\Big]. \label{1}%
\end{align}
Using the definition of the Hamiltonian $H$, we get
\begin{align}
I_{1}  &  =\mathbb{E}\Big[\int_{0}^{T}\int_{D}\Big(H_{0}(t,x)-\hat{H}%
_{0}(t,x)\Big)dxdt\nonumber\label{I1}\\
&  -\int_{0}^{T}\int_{D}\left\{  \widehat{p}(t,x)\widetilde{b}(t,x)+\widehat
{q}(t,x)\widetilde{\sigma}(t,x)\right\}  dxdt\Big].\nonumber
\end{align}
Summing the above we end up with
\begin{align*}
&  J(\xi)-J(\hat{\xi})=I_{1}+I_{2}+I_{3}\\
&  \leq\mathbb{E}\Big[{\int_{0}^{T}\int_{D}}\Big(H_{0}(t,x)-\widehat{H}%
_{0}(t,x)-\tilde{u}(t,x)\left\{  \frac{\partial\widehat{H}_{0}}{\partial
u}(t,x)+\overline{\nabla}_{\varphi}^{\ast}\widehat{H}_{0}(t,x)\right\}
\Big)dxdt\\
&  +\Big(H_{1}(t,x)\xi(dt,x)-\widehat{H}_{1}(t,x)\widehat{\xi}(t,x)-\tilde
{u}(t,x)\left\{  \frac{\partial\widehat{H}_{1}}{\partial u}(t,x)+\overline
{\nabla}_{\varphi}^{\ast}\widehat{H}_{1}(t,x)\right\}  \widehat{\xi
}(dt,x)\Big)dx\Big]\\
&  \leq\mathbb{E}\Big[{\int_{0}^{T}\int_{D}}\left\langle \nabla_{\xi}%
\widehat{H}(t,x),\xi(dt,x)-\widehat{\xi}(dt,x)\right\rangle dx\Big].
\end{align*}
\label{sum}\newline By the maximum condition of $H$ (\ref{max-cond}), we have%
\[
J(\xi)-J(\widehat{\xi})\leq0.
\]
\fproof

\subsection{A necessary maximum principle}

The concavity conditions in the sufficient maximum principle imposed on the
involved coefficients are not always satisfied. Hence, we will derive now a
necessary optimality conditions which do not require such an assumptions. We
shall first need the following Lemmas:\newline For $\xi\in\mathcal{A}$, we
let\textbf{ }$\mathcal{V}(\xi)$ denote the set of adapted processes
$\zeta(dt,x)$ of finite variation with respect to $t$, such that there exists
$\delta=\delta(\xi)>0$, such that $\xi+y\zeta\in\mathcal{A}$ for all
$y\in\lbrack0,\delta].$

\begin{lemma}
\label{lem-der} Let $\xi(dt,x)\in\mathcal{A}$ and choose $\zeta(dt,x)\in
\mathcal{V}(\xi)$. Define the derivative process
\begin{equation}
\mathcal{Z}(t,x):=\lim_{\epsilon\rightarrow0^{+}}\frac{1}{\epsilon}%
(u^{\xi+\epsilon\zeta}(t,x)-u^{\xi}(t,x)). \label{z}%
\end{equation}
Then $\mathcal{Z}$ satisfies the following singular linear SPDE
\begin{equation}
\left\{
\begin{array}
[c]{l}%
d\mathcal{Z}(t,x)=A_{x}\mathcal{Z}(t,x)dt+\left(  \frac{\partial b}{\partial
u}(t,x)\mathcal{Z}(t,x)+\left\langle \nabla_{\varphi}b(t,x),\mathcal{Z}%
(t,\cdot)\right\rangle \right)  dt\\
\qquad\qquad+\left(  \frac{\partial\sigma}{\partial u}(t,x)\mathcal{Z}%
(t,x)+\left\langle \nabla_{\varphi}\sigma(t,x),\mathcal{Z}(t,\cdot
)\right\rangle \right)  dB(t)\\
\qquad\qquad+f(t,x)\zeta(dt,x)\;;\quad(t,x)\in\lbrack0,T]\times D,\\
\mathcal{Z}(t,x)=0;\quad(t,x)\in(0,T)\times\partial D,\\
\mathcal{Z}(0,x)=0\;;\quad x\in D.
\end{array}
\right.  \label{der-proc}%
\end{equation}

\end{lemma}

\begin{lemma}
\label{der-j} Let $\xi(dt,x)\in\mathcal{A}$ and $\zeta(dt,x)\in\mathcal{V}%
(\xi)$. Put $\eta=\xi+\epsilon\zeta;\epsilon\in\lbrack0,\delta(\xi)]$. Then
\begin{align*}
\lim_{\epsilon\rightarrow0^{+}}  &  \frac{1}{\epsilon}(J(\xi+\epsilon
\zeta)-J(\xi))\\
&  =\mathbb{E}\left[  \int_{0}^{T}\int_{D}\left\{  f(t,x)p(t,x)+h_{1}%
(t,x,u(t,x),u(t,\cdot))\right\}  d\zeta(dt,x))\right]  .
\end{align*}

\end{lemma}

\dproof
By (\ref{ju}) and (\ref{der-proc}), we have
\begin{align}
&  \lim_{\epsilon\rightarrow0^{+}}\frac{1}{\epsilon}(J(\xi+\epsilon
\zeta)-J(\xi))\nonumber\\
&  =\mathbb{E}\left[  \int_{0}^{T}\int_{D}\left\{  \frac{\partial h_{0}%
}{\partial u}(t,x)\mathcal{Z}(t,x)+\left\langle \nabla_{\varphi}%
h_{0}(t,x),\mathcal{Z}(t,\cdot)\right\rangle \right\}  dxdt\right. \nonumber\\
&  +\int_{D}\left\{  \frac{\partial g}{\partial u}(T,x)\mathcal{Z}%
(T,x)+\left\langle \nabla_{\varphi}g(T,x),\mathcal{Z}(T,\cdot)\right\rangle
\right\}  dx\nonumber\\
&  +\int_{0}^{T}\int_{D}\left\{  \frac{\partial h_{1}}{\partial u}%
(t,x)\mathcal{Z}(t,x)+\left\langle \nabla_{\varphi}h_{1}(T,x),\mathcal{Z}%
(T,\cdot)\right\rangle \right\}  d\xi(t,x)\nonumber\\
&  \left.  +\int_{0}^{T}\int_{D}h_{1}(t,x)d\zeta(t,x)\right] \nonumber\\
&  =\mathbb{E}\left[  \int_{0}^{T}\int_{D}\left\{  \frac{\partial h_{0}%
}{\partial u}(t,x)\mathcal{Z}(t,x)+\overline{\nabla}_{\varphi}^{\ast}%
h_{0}(t,x)\mathcal{Z}(t,x)\right\}  dxdt\right. \nonumber\\
&  +\int_{D}\left\{  \frac{\partial g}{\partial u}(T,x)\mathcal{Z}%
(T,x)+\overline{\nabla}_{\varphi}^{\ast}g(T,x)\mathcal{Z}(T,x)\right\}
dxdt\nonumber\\
&  +\int_{0}^{T}\int_{D}\left\{  \frac{\partial h_{1}}{\partial u}%
(t,x)\mathcal{Z}(t,x)+\overline{\nabla}_{\varphi}^{\ast}h_{1}(T,x)\mathcal{Z}%
(T,x)\right\}  \xi(dt,x)dx\nonumber\\
&  \left.  +\int_{0}^{T}\int_{D}h_{1}(t,x)\zeta(dt,x)dx\right]  . \label{j}%
\end{align}
Using the definition \eqref{Ham} of the Hamiltonian, yields%
\begin{align}
&  \mathbb{E}\left[  \int_{0}^{T}\int_{D}\left\{  \frac{\partial h_{0}%
}{\partial u}(t,x)\mathcal{Z}(t,x)+\overline{\nabla}_{\varphi}^{\ast}%
h_{0}(t,x)\mathcal{Z}(t,x)\right\}  dxdt\right] \nonumber\\
&  =\mathbb{E}\left[  \int_{D}\int_{0}^{T}\left\{  \frac{\partial H_{0}%
}{\partial u}(t,x)\mathcal{Z}(t,x)+\overline{\nabla}_{\varphi}^{\ast}%
H_{0}(t,x)\mathcal{Z}(t,x)\right\}  dxdt\right. \nonumber\\
&  -\int_{0}^{T}\int_{D}\left\{  p(t,x)\left(  \frac{\partial b}{\partial
u}(t,x)\mathcal{Z}(t,x)+\overline{\nabla}_{\varphi}^{\ast}b(t,x)\mathcal{Z}%
(t,x)\right)  \right. \nonumber\\
&  \left.  +q(t,x)\left(  \frac{\partial\sigma}{\partial u}(t,x)\mathcal{Z}%
(t,x)+\overline{\nabla}_{\varphi}^{\ast}\sigma(t,x)\mathcal{Z}(t,x)\right)
\right\}  dxdt, \label{H}%
\end{align}
where we have used the simplified notation
\[
\frac{\partial H}{\partial u}(t,x)=\frac{\partial H}{\partial u}%
(t,x,u(t,x),u(t,\cdot),p(t,x),q(t,x))
\]
etc.\newline Applying the Itô formula to $p(T,x)\mathcal{Z}(T,x)$, we get
\begin{align}
&  \mathbb{E}\Big[\int_{D}\Big\{\frac{\partial g}{\partial u}(T,x)\mathcal{Z}%
(T,x)+\left\langle \nabla_{\varphi}g(T,x),\mathcal{Z}(T,\cdot)\right\rangle
\Big\}dx\Big]\nonumber\\
&  =\mathbb{E}\Big[\int_{D}p(T,x)\mathcal{Z}(T,x)dx\Big]\nonumber\\
&  =\mathbb{E}\Big[\int_{0}^{T}\int_{D}\Big\{p(t,x)\Big(A_{x}\mathcal{Z}%
(t,x)+\frac{\partial b}{\partial u}(t,x)\mathcal{Z}(t,x)+\left\langle
\nabla_{\varphi}b(t,x),\mathcal{Z}(t,\cdot)\right\rangle \Big)\nonumber\\
&  -A_{x}^{\ast}p(t,x)\mathcal{Z}(t,x)+\Big(\frac{\partial\sigma}{\partial
u}(t,x)\mathcal{Z}(t,x)+\left\langle \nabla_{\varphi}\sigma(t,x),\mathcal{Z}%
(t,\cdot)\right\rangle \Big)q(t,x)\Big\}dxdt\nonumber\\
&  \quad+\int_{0}^{T}\int_{D}f(t,x)p(t,x)\zeta(dt,x)dx\nonumber\\
&  -\int_{D}\int_{0}^{T}\Big(\frac{\partial H_{0}}{\partial u}(t,x)+\overline
{\nabla}_{\varphi}^{\ast}H_{0}(t,x)\Big)\mathcal{Z}(t,x)dxdt\nonumber\\
&  -\int_{D}\int_{0}^{T}\Big(\frac{\partial H_{1}}{\partial u}(t,x)+\overline
{\nabla}_{\varphi}^{\ast}H_{1}(t,x)\Big)\mathcal{Z}(t,x)\xi(dt,x)dx\Big].
\label{g}%
\end{align}
Since $p(t,x)=\mathcal{Z}(t,x)=0$ for $x\in\partial D$, we deduce that
\[
\int_{D}p(t,x)A_{x}\mathcal{Z}(t,x)dx=\int_{D}A_{x}^{\ast}p(t,x)\mathcal{Z}%
(t,x)dx.
\]
Therefore, substituting (\ref{g}) and (\ref{H}) into (\ref{j}), we get
\begin{align*}
\lim_{\epsilon\rightarrow0^{+}}  &  \frac{1}{\epsilon}(J(\xi+\epsilon
\zeta)-J(\xi))\\
&  =\mathbb{E}\left[  \int_{0}^{T}\int_{D}\left\{  f(t,x)p(t,x)+h_{1}%
(t,x)\right\}  \zeta(dt,x)dx\right]  .
\end{align*}
\fproof

We can now state our necessary maximum principle:

\begin{theorem}
[Necessary Maximum Principle]\label{ness}

(i) Suppose $\xi^{\ast}\in\mathcal{A}$ is optimal, i.e.
\begin{equation}
\max_{\xi\in\mathcal{A}}J(\xi)=J(\xi^{\ast}). \label{eq230}%
\end{equation}
Let $u^{\ast},(p^{\ast},q^{\ast})$ be the corresponding solution of
(\ref{eq2.1}) and (\ref{p}), respectively, and assume that \eqref{z} holds
with $\xi=\xi^{\ast}$. Then
\begin{equation}
f(t,x)p^{\ast}(t,x)+h_{1}(t,x,u^{\ast}(t,x),u^{\ast}(t,\cdot))\leq0\quad\text{
for all }t,x\in\lbrack0,T]\times D,\text{ a.s.,} \label{eq231}%
\end{equation}
and
\begin{equation}
\{f(t,x)p^{\ast}(t,x)+h_{1}(t,x,u^{\ast}(t,x),u^{\ast}(t,\cdot))\}\xi^{\ast
}(dt,x)=0\quad\text{ for all }t,x\in\lbrack0,T]\times D,\text{ a.s.}
\label{eq232}%
\end{equation}

(ii) Conversely, suppose that there exists $\hat{\xi}\in\mathcal{A},$ such
that the corresponding solutions $\widehat{u}(t,x),(\widehat{p}(t,x),\widehat
{q}(t,x))$ of (\ref{eq2.1}) and (\ref{p}), respectively, satisfy
\begin{equation}
f(t,x)\widehat{p}(t,x)+h_{1}(t,x,\widehat{u}(t,x),\widehat{u}(t,\cdot
))\leq0\quad\text{ for all }t,x\in\lbrack0,T]\times D,\text{ a.s.}
\label{eq233}%
\end{equation}
and
\begin{equation}
\left\{  f(t,x)\widehat{p}(t,x)+h_{1}(t,x,\widehat{u}(t,x),\widehat{u}%
(t,\cdot)\right\}  \widehat{\xi}(dt,x)=0\quad\text{ for all }t,x\in
\lbrack0,T]\times D,\text{ a.s.} \label{eq234}%
\end{equation}
Then $\widehat{\xi}$ is a directional sub-stationary point for $J(\cdot)$, in
the sense that
\begin{equation}
\lim_{\epsilon\rightarrow0^{+}}\frac{1}{\epsilon}\left(  J(\widehat{\xi
}+\epsilon\zeta)-J(\widehat{\xi})\right)  \leq0\quad\text{ for all }\zeta
\in\mathcal{V}(\hat{\xi}). \label{eq235}%
\end{equation}

\end{theorem}

\dproof
The proof is just a consequence of Lemma \ref{der-j}
and Theorem 3 in Øksendal \textit{et al} \cite{OSZ}.
\fproof

\section{Application to Optimal Harvesting}

We now return to the problem of optimal harvesting from a fish population in a
lake $D$ stated in the Introduction. Thus we suppose the density $u(t,x)$ of
the population at time $t\in\lbrack0,T]$ and at the point $x\in D$ is given by
the stochastic reaction-diffusion equation
\begin{equation}
\left\{
\begin{array}
[c]{l}%
du(t,x)=\left[  \dfrac{1}{2}\Delta u(t,x)+\alpha\bar{u}(t,x)\right]  dt+\beta
u(t,x)dB(t)-\lambda_{0}u(t,x)\xi(dt,x);\quad(t,x)\in(0,T)\times D,\\
\text{ }u(0,x)=u_{0}(x)>0;\quad x\in D,\\
\text{ }u(t,x)=u_{1}(t,x)\geq0;\quad(t,x)\in(0,T)\times\partial D,
\end{array}
\right.  \label{spde}%
\end{equation}
where $\lambda_{0}>0$ is a constant and, as in \eqref{G},
\[
\bar{u}(t,x)=\frac{1}{V(K_{\theta})}\int_{K_{\theta}}u(x+y)dy.
\]
The performance criterion is assumed to be
\[
J(\xi)=\mathbb{E}\left[  \int_{D}\int_{0}^{T}h_{10}(t,x)u(t,x)\xi
(dt,x)dx+\int_{D}g_{0}(T,x)u(T,x)dx\right]  ,
\]
where $h_{10}>0$ and $g_{0}>0$ are given deterministic functions. We can
interpret $\xi(dt,x)$ as the harvesting effort at $x$.

\begin{problem}
We want to find $\hat{\xi} \in\mathcal{A}$ such that $\sup_{\xi\in\mathcal{A}}
J(\xi) = J(\hat{\xi}).$ \label{problem}
\end{problem}

\noindent In this case the Hamiltonian is
\begin{align*}
&  H(t,x,u,\bar{u},p,q)(dt,\xi(dt,x))\\
&  =(\alpha\bar{u}p+\beta uq)dt+[-\lambda_{0}p+h_{10}(t,x)]u\xi(dt,x).
\end{align*}
Recall that for the map $L:L^{2}(D)\mapsto L^{2}(D)$ given by $L(u)=\bar{u}$
we know that
\[
\overline{\nabla}_{\varphi}^{\ast}L=\frac{V((x+K_{\theta})\cap D)}%
{V(K_{\theta})}.
\]
See Example 3.1 in Agram \textit{et al} \cite{AHO}. Therefore the adjoint
equation is
\begin{equation}
\left\{
\begin{array}
[c]{l}%
dp(t,x)=-\left[  \dfrac{1}{2}\Delta p(t,x)+\alpha p(t,x)\frac{V((x+K_{\theta
})\cap D)}{V(K_{\theta})}+\beta q(t,x)\right]  dt\\
\quad\quad\quad\quad+[\lambda_{0}-h_{10}(t,x)]\xi(dt,x)+q(t,x)dB(t,x);\quad
(t,x)\in(0,T)\times D,\\
p(T,x)=g_{0}(T,x);\quad x\in D,\\
p(t,x)=0;\quad(t,x)\in(0,T)\times\partial D.
\end{array}
\right.  \label{bspde}%
\end{equation}
The variational inequalities for an optimal control $\hat{\xi}(dt,x)$ and the
associated $\hat{p}$ are:
\begin{align*}
&  [-\lambda_{0}\hat{p}(t,x)+h_{10}(t,x)]\widehat{u}(t,x)\xi(dt,x)\\
&  \leq\lbrack-\lambda_{0}\hat{p}(t,x)+h_{10}(t,x)]\widehat{u}(t,x)\hat{\xi
}(dt,x);\quad(t,x)\in\lbrack0,T]\times D,\text{ for all }\xi.
\end{align*}
We claim that
\begin{equation}
u(t,x)>0\text{ for all }(t,x)\in\lbrack0,T]\times D. \label{eq3.5}%
\end{equation}
Suppose this claim is proved. Then, choosing first $\xi=2\hat{\xi}$ and then
$\xi=\frac{1}{2}\hat{\xi}$ in the above we obtain that
\[
\left[  \hat{p}(t,x)-\frac{1}{\lambda_{0}}h_{10}(t,x)\right]  \hat{\xi
}(dt,x)=0;\quad(t,x)\in\lbrack0,T]\times D.
\]
In addition we get that
\[
\left[  \hat{p}(t,x)-\frac{1}{\lambda_{0}}h_{10}(t,x)\right]  \hat{\xi
}(dt,x)\leq0;
\]
which implies that $\hat{p}(t,x)-\frac{1}{\lambda_{0}}h_{10}(t,x)\leq0$
always.\newline Summarising, we have proved the following:

\begin{theorem}
Suppose that $\widehat{u} > 0$ and $(\hat{p},\hat{\xi})$ satisfies the
following variational inequality
\begin{equation}
\label{eq3.8}\max\Big\{ \hat{p}(t,x)-\frac{1}{\lambda_{0}}h_{10}(t,x), -
\hat{\xi}(dt,x) \big \} = 0;\quad(t,x)\in\lbrack0,T]\times D.
\end{equation}
Then $\hat{\xi}$ is an optimal singular control for the space-mean SPDE
singular control problem \eqref{problem}
\end{theorem}

\noindent We see that this, together with \eqref{bspde} constitute a reflected
BSPDE, albeit of a slightly different type than the one that will be discussed
in the next section. \newline We summerize the above in the following:

\begin{theorem}
\begin{description}
\item[(a)] Suppose $\xi(dt,x)\in\mathcal{A}$ is an optimal singular control
for the harvesting problem
\[
\underset{\xi\in\mathcal{A}}{\sup}\mathbb{E}\left[  \int_{D}\int_{0}^{T}%
h_{1}(t,x)\xi(dt,x)dx+\int_{D}g_{0}(T,x)u(T,x)dx\right]  ,
\]
where $u(t,x)$ is given by the SPDE \eqref{spde}. Then $\xi(dt,x)$ solves the
reflected BSPDE \eqref{bspde}, \eqref{eq3.8}. \vskip0.1cm

\item[(b)] Conversely, suppose $(p,q,\xi)$ is a solution of the reflected
BSPDE \eqref{bspde}, \eqref{eq3.8}. Then $\xi(dt,x)$ is an optimal control for
the problem to maximize the performance \eqref{eq248}.
\end{description}
\end{theorem}

Heuristically we can interpret the optimal harvesting strategy as follows:

\begin{itemize}
\item As long as $p(t,x)<\frac{1}{\lambda_{0}}h_{1}(t,x)$, we do nothing.

\item If $p(t,x)=\frac{1}{\lambda_{0}}h_{1}(t,x)$, we harvest immediately from
$u(t,x)$ at a rate $\xi(dt,x)$ which is exactly enough to prevent $p(t,x)$
from dropping below $\frac{1}{\lambda_{0}}h_{1}(t,x)$ in the next moment.

\item If $p(t,x)>\frac{1}{\lambda_{0}}h_{1}(t,x)$, we harvest immediately what
is necessary to bring $p(t,x)$ up to the level of $\frac{1}{\lambda_{0}}%
h_{1}(t,x).$
\end{itemize}

\begin{remark}
Note that if $p(t,x)=\frac{1}{\lambda_{0}}h_{10}(t,x)$ and
\[
\lambda_{0}>h_{10}(t,x),
\]
then an immediate harvesting of an amount $\Delta\xi>0$ from $u(t,x)$ produces
an immediate decrease in the process $p(t,x)$ and hence pushes $p(t,x)$ below
$\frac{1}{\lambda_{0}}h_{10}(t,x).$ This follows from the comparison theorem
for reflected BSPDEs of the type \eqref{bspde}.
\end{remark}

\section{Existence and uniqueness of solutions of space-mean reflected
backward SPDEs}

Let $W,H$ be two separable Hilbert spaces such that $W$ is continuously,
densely imbedded in $H$. Identifying $H$ with its dual we have
\[
W\subset H\approxeq H^{\ast}\subset W^{\ast},
\]
where we have denoted by $W^{\ast}$ the topological dual of $V$. Let $A$ be a
bounded linear operator from $W$ to $W^{\ast}$ satisfying the following
Gårding inequality (coercivity hypothesis): There exist constants $\alpha>0$
and $\lambda\geq0$ so that
\begin{equation}
2\langle Au,u\rangle+\lambda||u||_{H}^{2}\geq\alpha||u||_{W}^{2};\,\,\,\text{
for all }u\in W, \label{COE}%
\end{equation}
where $\langle Au,u\rangle=Au(u)$ denotes the action of $Au\in W^{\ast}$ on
$u\in W$ and $||\cdot||_{H}$ (respectively $\Vert\cdot\Vert_{W}$) the norm
associated to the Hilbert space $H$ (respectively $W$). We will also use the
following spaces:

\begin{itemize}
\item $L^{2}(D)$ is the set of all Lebesgue measurable $Y:D\rightarrow%
%TCIMACRO{\U{211d} }%
%BeginExpansion
\mathbb{R}
%EndExpansion
,$ such that%
\[
||Y||_{L^{2}(D)}:=\left(  \int_{D}|Y(x)|^{2}dx\right)  ^{\frac{1}{2}}<\infty.
\]

\item $L^{2}(H)$ is the set of $\mathcal{F}_{T}$-measurable $H$-valued random
variables $\varsigma$ such that $\mathbb{E}[||\varsigma||_{H}^{2}]<\infty$.
\end{itemize}

\noindent We let $W:=W^{1,2}(D)$ and $H=L^{2}(D).$

\noindent Denote by $L(t,x)$ the barrier which is a measurable function that
is differentiable in time $t$ and twice differentiable in space $x,$ such
that
\[
\int_{0}^{T}\int_{D}L^{\prime}(t,x)^{2}dxdt<\infty,\text{ }\int_{0}^{T}%
\int_{D}|\Delta L(t,x)|^{2}dxdt<\infty.
\]

\noindent$\eta$ is a $H$-valued continuous process, nonnegative, nondecreasing in
$t$ and $\eta(0,x)=0.$

\noindent We now consider the adjoint equation (\ref{p}) as a reflected
backward stochastic evolution equation
\begin{equation}
\left\{
\begin{array}
[c]{l}%
dY(t,x)=-AY(t,x)dt-F(t,Y(t,x),\overline{Y}(t,x),Z(t,x),\overline{Z}(t,x))dt\\
\quad\quad\quad\quad+Z(t,x)dB(t)-\eta(dt,x),t\in(0,T),\\
Y(t,x)\geq L(t,x),\\
\int_{0}^{T}\int_{D}(Y(t,x)-L(t,x))\eta(dt,x)dx=0,\\
Y(T,x)=\phi(x);\quad\text{a.s.,}%
\end{array}
\right.  \label{r-BSPDE}%
\end{equation}
where $Y(t,x)$ stands for the $W$-valued continuous process $Y(t,x)$ and the
solution of equation $(\ref{r-BSPDE})$ is understood as an equation in the
dual space $W^{\ast}$ of $W$.\newline We mean by $dY(t,x)$ the differential
operator with respect to $t$, while $A_{x}$ is the partial differential
operator with respect to $x$, and \newline%
\begin{align*}
\overline{Y}(t,x)  &  =G(x,Y)=\frac{1}{V}\int_{K_{\theta}}Y(x+\rho)d\rho,\\
\overline{Z}(t,x)  &  =G(x,Z)=\frac{1}{V}\int_{K_{\theta}}Z(x+\rho)d\rho.
\end{align*}
The following result is essential due to Agram \textit{et al }\cite{AHO}:

\begin{lemma}
\label{lemexi}For all $\varphi\in H$ we have
\begin{equation}
||G(\cdot,\varphi)||_{H}\leq||\varphi||_{H}. \label{eq1.3}%
\end{equation}

\end{lemma}

\noindent We shall now state and prove our main result of existence and
uniqueness of solutions to reflected BSPDE.

\begin{theorem}
[Existence and uniqueness of solutions]The space-mean reflected BSPDE
$(\ref{r-BSPDE})$ has a unique solution $(Y(t,x),Z(t,x),\eta(t,x))\in W\times
L^{2}(D,\mathbb{R}^{m})\times H$\textit{-valued progressively measurable
process, }provided that the following assumptions hold:

\begin{description}
\item[(i)] The terminal condition $\phi$ is $\mathcal{F}_{T}$-measurable
random variable and satisfies%
\[
\mathbb{E}\left[  ||\phi||_{H}^{2}\right]  <\infty.
\]

\item[(ii)] There exists a constant $C>0$ such that
\begin{align*}
&  ||F(t,y_{1},\overline{y}_{1},z_{1},\overline{z}_{1})-F(t,y_{2},\overline
{y}_{2},z_{2},\overline{z}_{2})||_{H}\\
&  \leq C\left(  ||y_{1}-y_{2}||_{H}+||\overline{y}_{1}-\overline{y}_{2}%
||_{H}+||z_{1}-z_{2}||_{H}+||\overline{z}_{1}-\overline{z}_{2}||_{H}\right)  ,
\end{align*}

\end{description}

for all $t,y_{i},\overline{y}_{i},z_{i},\overline{z}_{i};i=1,2.$
\end{theorem}

\dproof
For the proof of the theorem, we introduce the
penalized backward SPDEs:
\begin{equation}
\left\{
\begin{array}
[c]{l}%
dY^{n}(t,x)=-AY^{n}(t)dt-F(t,Y^{n}(t,x),\overline{Y}^{n}(t,x),Z^{n}%
(t,x),\overline{Z}^{n}(t,x))dt\\
\quad\quad\quad\quad+Z^{n}(t,x)dB(t)-n(Y^{n}(t,x)-L(t,x))^{-}dt,\quad
t\in(0,T),\\
Y^{n}(T,x)=\phi(x)\quad\text{a.s.}%
\end{array}
\right.  \label{2.2}%
\end{equation}
According to Agram \textit{et al} \cite{AHO}, the solution $(Y^{n},Z^{n})$ of
the above equation (\ref{2.2}) exists and is unique. We are going to show that
$(Y^{n},Z^{n})_{n\geq1}$ forms a Cauchy sequence, i.e.,%
\[
\lim_{n,m\rightarrow\infty}\mathbb{E}\left[  \sup_{0\leq t\leq T}%
|Y^{n}(t)-Y^{m}(t)|_{H}^{2}\right]  =0,
\]%
\[
\lim_{n,m\rightarrow\infty}\mathbb{E}\left[  \int_{0}^{T}||Y^{n}%
(t)-Y^{m}(t)||_{W}^{2}dt\right]  =0,
\]%
\[
\lim_{n,m\rightarrow\infty}\mathbb{E}\left[  \int_{0}^{T}|Z^{n}(t)-Z^{m}%
(t)|_{L^{2}(D,\mathbb{R}^{m})}^{2}dt\right]  =0.
\]
Applying Itô's formula, it follows that%

\begin{align*}
&  |Y^{n}(t)-Y^{m}(t)|_{H}^{2}\\
&  =2\int_{t}^{T}\left\langle Y^{n}(s)-Y^{m}(s),A(Y^{n}(s)-Y^{m}%
(s))\right\rangle ds\\
&  +2\int_{t}^{T}\left\langle Y^{n}(s)-Y^{m}(s),\right. \\
&  \left.  F(s,Y^{n}(s),\overline{Y}^{n}(s),Z^{n}(s),\overline{Z}%
^{n}(s))-F(s,Y^{m}(s),\overline{Y}^{m}(s),Z^{m}(s)),\overline{Z}%
^{m}(s))\right\rangle ds\\
&  -2\int_{t}^{T}\left\langle Y^{n}(s)-Y^{m}(s),Z^{n}(s)-Z^{m}(s)\right\rangle
dB(s)\\
&  +2\int_{t}^{T}\left\langle Y^{n}(s)-Y^{m}(s),n(u^{n}(s)-L(s))^{-}%
-m(Y^{m}(s)-L(s))^{-}\right\rangle ds\\
&  -\int_{t}^{T}|Z^{n}(s)-Z^{m}(s)|_{L^{2}(D,\mathbb{R}^{m})}^{2}ds.
\end{align*}
Now we estimate each of the terms on the right side:%

\begin{align}
&  2\int_{t}^{T}\left\langle Y^{n}(s)-Y^{m}(s),A(Y^{n}(s)-Y^{m}%
(s))\right\rangle ds\nonumber\\
&  \leq\lambda\int_{t}^{T}||Y^{n}(s)-Y^{m}(s)||_{H}^{2}ds-\alpha\int_{t}%
^{T}||Y^{n}(s)-Y^{m}(s)||_{V}^{2}ds. \label{4.5}%
\end{align}
By the Lipschitz continuity of $b$ and the inequality $ab\leq\varepsilon
a^{2}+C_{\varepsilon}b^{2}$, together with inequality (\ref{eq1.3}), one has
\begin{align}
&  2\int_{t}^{T}\left\langle Y^{n}(s)-Y^{m}(s),\right. \nonumber\\
&  \left.  F(s,Y^{n}(s),\overline{Y}^{n}(s),Z^{n}(s),\overline{Z}%
^{n}(s))-F(s,Y^{m}(s),\overline{Y}^{m}(s),Z^{m}(s)),\overline{Z}%
^{m}(s))\right\rangle ds\nonumber\\
&  \leq C\int_{t}^{T}|Y^{n}(s)-Y^{m}(s)|_{H}^{2}ds+\frac{1}{2}\int_{t}%
^{T}|Z^{n}(s)-Z^{m}(s)|_{L^{2}(D,\mathbb{R}^{m})}^{2}ds. \label{4.6}%
\end{align}
It follows from (\ref{4.5}) and (\ref{4.6}) that
\begin{align*}
&  \mathbb{E}[|Y^{n}(t)-Y^{m}(t)|_{W}^{2}]+\frac{1}{2}\mathbb{E}\left[
\int_{t}^{T}|Z^{n}(s)-Z^{m}(s)|_{L^{2}(D,\mathbb{R}^{m})}^{2}ds\right] \\
&  +\mathbb{E}\left[  \int_{t}^{T}||Y^{n}(s)-Y^{m}(s)||_{W}^{2}ds\right] \\
&  \leq C\int_{t}^{T}\mathbb{E}[|Y^{n}(s)-Y^{m}(s)|_{W}^{2}]ds+C^{\prime
}\left(  \frac{1}{n}+\frac{1}{m}\right)  .
\end{align*}
Gronwall inequality, yields
\begin{equation}
\lim_{n,m\rightarrow\infty}\left\{  \mathbb{E}[|Y^{n}(t)-Y^{m}(t)|_{H}%
^{2}]+\frac{1}{2}\mathbb{E}\left[  \int_{t}^{T}|Z^{n}(s)-Z^{m}(s)|_{L^{2}%
(D,\mathbb{R}^{m})}^{2}ds\right]  \right\}  =0, \label{4.7}%
\end{equation}
and%
\[
\lim_{n,m\rightarrow\infty}\mathbb{E}\left[  \int_{t}^{T}||Y^{n}%
(s)-Y^{m}(s)||_{H}^{2}ds\right]  =0.
\]
By inequality (\ref{4.7}) and the Burkholder inequality we get
\[
\lim_{n,m\rightarrow\infty}\mathbb{E}\left[  \sup_{0\leq t\leq T}%
|Y^{n}(t)-Y^{m}(t)|_{H}^{2}\right]  =0.
\]

\noindent Under the conditions of Theorem 4.2 and by Lemma 5 in Øksendal \textit{et al} \cite{OSZ}, there exists
a constant $C,$ such that
\begin{equation}
\mathbb{E}\left[  \int_{0}^{T}\int_{D}((Y^{n}(t,x)-L(t,x))^{-})^{2}%
dxdt\right]  \leq\frac{C}{n^{2}}. \label{dp}%
\end{equation}

\noindent Denote by $Y(t,x)$, $Z(t,x)$ the limit of $Y^{n}$ and $Z^{n}$,
respectively. Put
\[
\overline{\eta}^{n}(t,x)=n(Y^{n}(t,x)-L(t,x))^{-}.
\]
Inequality (\ref{dp}) implies that $\overline{\eta}^{n}(t,x)$ admits a
non-negative weak limit, denoted by $\overline{\eta}(t,x)$, in the following
Hilbert space:
\[
\overline{H}=\left\{  h;\quad h\text{ is a }H\text{-valued adapted process,
such that}\quad\mathbb{E}\left[  \int_{0}^{T}|h(s)|_{H}^{2}ds\right]
<\infty\right\}  ,
\]
with inner product
\[
\left\langle h_{1},h_{2}\right\rangle _{\overline{H}}=\mathbb{E}\left[
\int_{0}^{T}\int_{D}h_{1}(t,x)h_{2}(t,x)dtdx\right]  .
\]
Set $\eta(t,x)=\int_{0}^{t}\overline{\eta}(s,x)ds$. Then $\eta$ is a
continuous $H$-valued process which is increasing in $t$. Letting
$n\rightarrow\infty$ in (\ref{2.2}) we obtain
\begin{align}
&  Y(t,x)\nonumber\\
&  =\phi(x)+\int_{t}^{T}AY(s,x)ds+\int_{t}^{T}F(s,Y(s,x),\overline
{Y}(s,x),Z(s,x),\overline{Z}(s,x))ds\nonumber\\
&  -\int_{t}^{T}Z(s,x)dB(s)+\eta(T,x)-\eta(t,x);\quad0\leq t\leq T.
\label{2.48}%
\end{align}
Inequality (\ref{dp}) and the Fatou Lemma imply that $\mathbb{E}\left[
\int_{t}^{T}\int_{D}((Y(s,x)-L(s,x))^{-})^{2}dxds\right]  =0$. In view of the
continuity of $Y$ in $t$, we conclude $Y(t,x)\geq L(t,x)$ a.e. in $x$, for
every $t\geq0$. Combining the strong convergence of $Y^{n}$ and the weak
convergence of $\bar{\eta}^{n}$, we also have
\begin{align}
&  \mathbb{E}\left[  \int_{0}^{T}\int_{D}(Y(s,x)-L(s,x))\eta(dt,x)dx\right]
\nonumber\label{2.49}\\
&  =\mathbb{E}\left[  \int_{0}^{T}\int_{D}(Y(s,x)-L(s,x))\overline{\eta
}(t,x)dtdx\right] \nonumber\\
&  \leq\lim_{n\rightarrow\infty}\mathbb{E}\left[  \int_{0}^{T}\int_{D}%
(Y^{n}(s,x)-L(s,x))\overline{\eta}^{n}(t,x)dtdx\right]  \leq0.
\end{align}
Hence,
\[
\int_{0}^{T}\int_{D}(Y(s,x)-L(s,x))\eta(dt,x)dx=0,\quad\text{a.s.}%
\]
We have shown that $(Y,Z,\eta)$ is a solution to the reflected backward SPDE
$(\ref{r-BSPDE})$. \vskip0.4cm

\noindent\textbf{Uniqueness}. Let $(Y_{1},Z_{1},\eta_{1})$, $(Y_{2},Z_{2}%
,\eta_{2})$ be two such solutions to equation $(\ref{r-BSPDE})$. By Itô's
formula, we have%

\begin{align}
&  |Y_{1}(t)-Y_{2}(t)|_{H}^{2}\nonumber\\
&  =2\int_{t}^{T}\left\langle Y_{1}(s)-Y_{2}(s),\Delta(Y_{1}(s)-Y_{2}%
(s))\right\rangle ds\nonumber\\
&  +2\int_{t}^{T}\left\langle Y_{1}(s)-Y_{2}(s),\right. \nonumber\\
&  \left.  F(s,Y_{1}(s),\overline{Y}_{1}(s),Z_{1}(s),\overline{Z}%
_{1}(s))-F(s,Y_{2}(s),\overline{Y}_{2}(s),Z_{2}(s),\overline{Z}_{2}%
(s))\right\rangle ds\nonumber\\
&  -2\int_{t}^{T}\left\langle Y_{1}(s)-Y_{2}(s),Z_{1}(s)-Z_{2}(s)\right\rangle
dB(s)\nonumber\\
&  +2\int_{t}^{T}\left\langle Y_{1}(s)-Y_{2}(s),\eta_{1}(ds)-\eta
_{2}(ds)\right\rangle \nonumber\\
&  -\int_{t}^{T}|Z_{1}(s)-Z_{2}(s)|_{L^{2}(D,\mathbb{R}^{m})}^{2}ds.
\label{2.50}%
\end{align}
Similar to the proof of existence, we have
\begin{equation}
2\int_{t}^{T}\left\langle Y_{1}(s)-Y_{2}(s),A(Y_{1}(s)-Y_{2}(s))\right\rangle
ds\leq0, \label{2.51}%
\end{equation}
and
\begin{align}
&  2\int_{t}^{T}\left\langle Y_{1}(s)-Y_{2}(s),\right. \nonumber\\
&  \left.  F(s,Y_{1}(s),\overline{Y}_{1}(s),Z_{1}(s),\overline{Z}%
_{1}(s))-F(s,Y_{2}(s),\overline{Y}_{2}(s),Z_{2}(s),\overline{Z}_{2}%
(s))\right\rangle ds\nonumber\\
&  \leq C\int_{t}^{T}|Y_{1}(s)-Y_{2}(s)|_{H}^{2}ds+\frac{1}{2}\int_{t}%
^{T}|Z_{1}(s)-Z_{2}(s)|_{L^{2}(D,\mathbb{R}^{m})}^{2}ds \label{2.52}%
\end{align}
On the other hand,
\begin{align}
&  2\mathbb{E}\left[  \int_{t}^{T}\left\langle Y_{1}(s)-Y_{2}(s),\eta
_{1}(ds)-\eta_{2}(ds)\right\rangle \right] \nonumber\\
&  =2\mathbb{E}\left[  \int_{t}^{T}\int_{D}(Y_{1}(s,x)-L(s,x))\eta
_{1}(ds,x)dx\right] \nonumber\\
&  -2\mathbb{E}\left[  \int_{t}^{T}\int_{D}(Y_{1}(s,x)-L(s,x))\eta
_{2}(ds,x)dx\right] \nonumber\\
&  +2\mathbb{E}\left[  \int_{t}^{T}\int_{D}(Y_{2}(s,x)-L(s,x))\eta
_{2}(ds,x)dx\right] \nonumber\\
&  -2\mathbb{E}\left[  \int_{t}^{T}\int_{D}(Y_{2}(s,x)-L(s,x))\eta
_{1}(ds,x)dx\right] \nonumber\\
&  \leq0. \label{2.53}%
\end{align}
Combining (\ref{2.50})-(\ref{2.53}) we arrive at
\begin{align*}
&  \mathbb{E}[|Y_{1}(t)-Y_{2}(t)|_{H}^{2}]+\frac{1}{2}\mathbb{E}\left[
\int_{t}^{T}|Z_{1}(s)-Z_{2}(s)|_{L^{2}(D,\mathbb{R}^{m})}^{2}ds\right] \\
&  \leq C\int_{t}^{T}\mathbb{E}[|Y_{1}(s)-Y_{2}(s)|_{H}^{2}]ds.
\end{align*}
Appealing to the Gronwall inequality, this implies
\[
Y_{1}=Y_{2},\quad Z_{1}=Z_{2}%
\]
which further gives $\eta_{1}=\eta_{2}$ from the equation they satisfy.
\fproof

%%%%%%%%%%%%%%%%%%%%%%%%%%%%%
%%%%%%%%%%%%%%%%%%%%%%%%%%%%%

\end{document}